\documentclass[11pt,twoside]{article}

\usepackage{epsfig,amsfonts,color}
\usepackage{amsmath}

\usepackage{amssymb, palatino, geometry,url}
\usepackage{algorithmic}
\usepackage[noresetcount,lined,boxed]{algorithm2e} 

\usepackage[colorlinks=true,linkcolor=blue,citecolor=blue,urlcolor=blue]{hyperref}
\usepackage{subcaption}
\usepackage{amsthm}

\usepackage{fancyhdr}
\newcommand{\btheta}{\boldsymbol{\theta}}

\newcommand{\amin}{\mathop{\mbox{argmin}}}

\newcommand{\be}{\begin{equation}}
\newcommand{\ee}{\end{equation}}
\newcommand{\bea}{\begin{eqnarray}}
\newcommand{\eea}{\end{eqnarray}}

\newcommand{\bvec}{\left(\begin{array}{c}}
	\newcommand{\evec}{\end{array}\right)}
\newcommand{\bsub}{\begin{subequations}}
	\newcommand{\esub}{\end{subequations}}

\usepackage{float}
\usepackage{bbm}
\usepackage[section]{placeins}
\usepackage{multirow}

\title{A Scalable Stochastic Programming Approach\\ for the Design of Flexible Systems}
\author{Joshua L. Pulsipher and Victor M. Zavala\thanks{Corresponding Author: victor.zavala@wisc.edu}\\
	{\small Department of Chemical and Biological Engineering}\\
	{\small \;University of Wisconsin, 1415 Engineering Dr, Madison, WI 53706, USA}}
\date{}

\begin{document}
\maketitle

\begin{abstract}
	We study the problem of designing systems in order to minimize cost while meeting a given flexibility target. Flexibility is attained by enforcing a joint chance constraint, which ensures that the system will exhibit feasible operation with a given target probability level. Unfortunately, joint chance constraints are complicated mathematical objects that often need to be reformulated using mixed-integer programming (MIP) techniques. In this work, we cast the design problem as a conflict resolution problem that seeks to minimize cost while maximizing flexibility. We propose a purely continuous relaxation of this problem that provides a significantly more scalable approach relative to MIP methods and show that the formulation delivers solutions that closely approximate the Pareto set of the original joint chance-constrained problem.  
\end{abstract}


\section{Problem Definition and Setting}

Flexibility is the ability of a system to maintain feasible operation in the face of externalities. A number of approaches have been proposed to quantify and analyze system flexibility and to design flexible systems \cite{grossmann2014evolution}. The flexibility index problem, first proposed in \cite{grossmann1983optimization}, for a given design $\mathbf{d}\in \mathbb{R}^{n_d}$ seeks to identify the largest uncertainty set $T(\delta)$ (where $\delta \in \mathbb{R}_+$ is a parameter that scales $T$) for which the system remains feasible. In other words, it seeks to find the largest uncertainty set under which there exists recourse system variables $\mathbf{z} \in \mathbb{R}^{n_z}$ that counteract the uncertain parameters $\btheta\in T(\delta)\subseteq \mathbb{R}^{n_\theta}$ in order to satisfy the system constraints $f_j(\mathbf{d}, \mathbf{z}, \btheta)\leq 0,\,j\in J$, and $h_i(\mathbf{d}, \mathbf{z}, \btheta)= 0,\,i\in I$. The flexibility index $F(\mathbf{d})$ is defined as:
\begin{equation}
\begin{aligned}
& F(\mathbf{d}) := & \max_{\delta\in\mathbb{R}_+} &&&\delta \\
&&\text{s.t.} &&& \max_{\btheta \in T(\delta)} \psi(\mathbf{d}, \btheta) \leq 0.
\end{aligned}
\label{eq:flexibility_index_definition}
\end{equation}
\noindent where $\mathbf{d} \in \mathbb{R}^{n_d}$ are design variables and $\psi(\mathbf{d}, \btheta)$ is a function that tests the feasibility of the system at a particular instance of $\btheta$. The feasibility function $\psi(\mathbf{d}, \btheta)$ is defined:

\begin{equation}
\begin{aligned}
& \psi(\mathbf{d}, \boldsymbol{\theta}) := & \min_{\mathbf{z}} & \max_{j \in J} \ f_j(\mathbf{d}, \mathbf{z}, \boldsymbol{\theta}) \\
&& \text{s.t.} & \ \ h_i(\mathbf{d}, \mathbf{z}, \boldsymbol{\theta}) = 0, & i \in I.
\end{aligned}
\label{eq:feasibility_function}
\end{equation}

\noindent Here, the system exhibits feasible operation at a particular instance $\btheta$ if $\psi(\mathbf{d}, \boldsymbol{\theta}) \leq 0$. Notably, the feasibility function is a minimax problem and thus can be cast as a standard optimization by the use of an upper bounding variable $u \in \mathbb{R}$:
\begin{equation}
\begin{aligned}
& \psi(\mathbf{d}, \boldsymbol{\theta}) := & \min_{\mathbf{z}} &&& u \\
&& \text{s.t.} &&& f_j(\mathbf{d}, \mathbf{z}, \boldsymbol{\theta})  \leq y, & j \in J \\
&&&&& h_i(\mathbf{d}, \mathbf{z}, \boldsymbol{\theta}) = 0, & i \in I.
\end{aligned}
\label{eq:feasibility_function_alt}
\end{equation}

Hence, the flexibility index problem simply seeks to determine the largest uncertainty space $T(\delta)$ over which feasible operation can be encountered. The $F$ index is a flexibility metric deterministic in nature that shares some common ground with robust optimization methods, Zhang et. al. discuss this relationship in detail for linear systems \cite{zhang2016relation}. Swaney and Grossmann proved in \cite{swaney1985index1} that this is equivalent to searching for the minimum $\delta$ along the boundary of the feasible region, provided that $T(\delta)$ is compact and the constraints $f_j(\mathbf{d}, \mathbf{z}, \mathbf{x}, \btheta)$ and $h_i(\mathbf{d}, \mathbf{z}, \mathbf{x}, \btheta)$ are continuous in $\mathbf{z}$, $\mathbf{x}$, and $\btheta$. Thus, the $F$ index can be expressed as a bi-level optimization problem under these conditions. Grossmann and Floudas leverage this observation to propose an active constraint mixed-integer program (MIP) that casts the inner problem in terms of its first-order Karush-Kuhn-Tucker (KKT) conditions  \cite{grossmann1987active}.

The flexibility index problem has given rise to the most common class of design problems in flexibility analysis where design variables $\mathbf{d}$ are selected to minimize cost and either ensure the feasibility at a fixed $F$ index or maximize the flexibility index problem \cite{grossmann1983optimization,grossmann1979optimum,pistikopoulos1988optimal,varvarezos1994multiperiod}. However, studies reported in the literature find the latter problem to be largely intractable due to its multi-objective complexity \cite{grossmann2014evolution, pistikopoulos1989optimal}. Thus, the design problem has traditionally been formulated such that it enforces the feasibility of an uncertainty set $T(\hat{F})$ fixed at a particular flexibility index $\hat{F}$. This optimal design problem is formalized 
\begin{equation}
\begin{aligned}
&&\amin_{\mathbf{d}, \btheta} &&& c(\mathbf{d}) \\
&&\text{s.t.} &&& \max_{\btheta \in T(\hat{F})} \psi(\mathbf{d}, \btheta) \leq 0\\
&&&&& \mathbf{d} \in D
\end{aligned}
\label{eq:set_design}
\end{equation}
where $c(\mathbf{d})$ is the cost function and $D$ denotes the set of all feasible $\mathbf{d}$ values. This is a straightforward approach for selecting $\mathbf{d}$, however the use of uncertainty sets conservatively models the flexibility of a given system design \cite{pulsipher2018MICP}. This means that solutions to the design problem with the use of uncertainty sets might be overly conservative and thus choose a greater cost than is necessary to achieve a certain degree of flexibility.

Straub and Grossmann proposed a probabilistic flexibility measure that they called the stochastic flexibility ($SF$) index. This index quantifies the probability of finding feasible operation \cite{straub1990integrated}. Here, the uncertain parameters $\btheta$ are modeled as random variables with associated probability density function $p:\mathbb{R}^{n_\theta}\to \mathbb{R}$. The stochastic flexibility index is defined:

\begin{equation}
SF(\mathbf{d}) := \int_{\boldsymbol{\theta}\in\Theta(\mathbf{d})} p(\boldsymbol{\theta}) d\boldsymbol{\theta}
\label{eq:SF_definition}
\end{equation}

\noindent where $\Theta (\mathbf{d}):=\{\boldsymbol{\theta}\,:\,\psi(\mathbf{d}, \boldsymbol{\theta}) \leq 0\}$ is the feasible set (projected onto the space of the uncertain parameters). Pistikopoulos and Mazzuchi proposed a similar definition in \cite{pistikopoulos1990novel}.  In \cite{pulsipher2018MICP}  it is noted that the SF index can also be expressed as:

\begin{align}
SF(\mathbf{d}) &= \mathbb{P}\left( \psi(\mathbf{d}, \boldsymbol{\theta}) \leq 0 \right) \nonumber\\
&= \mathbb{P}\left( \exists \ \mathbf{z} :  f_j(\mathbf{d}, \mathbf{z}, \boldsymbol{\theta}) \leq 0, \  h_i(\mathbf{d}, \mathbf{z}, \boldsymbol{\theta}) = 0,\, j \in J, i \in I \right).
\label{eq:SF_definition2}
\end{align}

The $SF$ index can be computed rigorously by evaluating  \eqref{eq:SF_definition} or \eqref{eq:SF_definition2} via Monte Carlo (MC) sampling. Here, the feasibility of each realization of $\btheta$ is assessed using the feasibility function $\psi(\mathbf{d},\btheta)$. Such an approach converges asymptotically but can require a large number of samples \cite{shapiro2013sample, robert2013monte}. This is a common drawback of sample average approaximation approaches for stochastic programming but is especially important in the context of flexibility analysis because limiting behavior is often found near the boundary of $\Theta (\mathbf{d})$. Interestingly, high-dimensional multivariate Gaussian random variables have been shown to contain most probability mass in a thin ellipsoidal shell from which the MC samples are generated \cite{scott2015multivariate}. Consequently, in some applications a moderate number of samples might be needed. Important sampling and sparse grid techniques have also been proposed recently to cover the uncertainty space more effectively and with this reduce the number of samples \cite{barrera2016chance,renteria2018optimal}.

The $SF$ index can be used as a metric to guide the design of flexible systems \cite{straub1993design}. Typically, one aims to find a design $\mathbf{d}$ that minimizes a design cost while ensuring that the system remains feasible with a given probability level $\alpha\in [0,1]$ (typically close to one): 
\begin{align} \label{eq:design}
\min_{\mathbf{d}\in \mathcal{D}} &\; c(\mathbf{d})\nonumber \\ 
\textrm{s.t.} &\; SF(\mathbf{d})\geq \alpha. 
\end{align}
here, $c:\mathbb{R}^{n_d}\to \mathbb{R}$ is the cost design function and $\mathcal{D}\subseteq\mathbb{R}^{n_d}$ is the design space. Note that this is an optimization problem with {\em joint chance constraints} that is computationally challenging to solve. In particular, the joint chance constraint often needs to be reformulated by using binary variables \cite{luedtke2008sample}. In this work, we provide a more scalable approach to address the design problem. The approach relies on the observation that the above problem provides a Pareto solution (for a specific value of $\alpha$) for the conflict resolution (multi-objective) problem: 
\begin{align} 
\min_{\mathbf{d}\in \mathcal{D}} \left\{c(\mathbf{d}),-SF(\mathbf{d})\right\}.
\end{align}
We will demonstrate that we can recover the Pareto set for this problem to high accuracy (which includes solutions of the design problem \eqref{eq:design}) by solving a {\em continuous} formulation. This thus provides a scalable approach to solve large-scale design problems. 

Section \ref{sec:derivation} outlines the derivation of the standard approaches and the proposed approach to address the design problem. Section \ref{sec:cases} provides analysis in conjunction with large case studies to demonstrate the properties of the formulations.

\section{Solution Approaches for the Design Problem} \label{sec:derivation}

This section outlines standard approaches to solve the conflict resolution problem using sample average approximations and mixed-integer reformulations. Here, we also discuss the continuous formulation used to approximate the solution of the design problem.

\subsection{Sample Average Approximation of the Stochastic Flexibility Index}

\noindent The $SF$ index (shown in \eqref{eq:SF_definition} and \eqref{eq:SF_definition2}) involves a high-dimensional integral and a projection into the feasible space of the system. This index can be  approximated via Monte Carlo (MC) sampling as:

\begin{equation}
SF(\mathbf{d}) = \mathbb{E}[\mathbbm{1}_{\boldsymbol{\theta} \in \Theta (\mathbf{d})}] \approx SF_K(\mathbf{d}):=\frac{1}{|K|} \sum_{k \in K} \mathbbm{1}_{\boldsymbol{\theta}^k \in \Theta (\mathbf{d})}
\label{eq:mc_approx}
\end{equation}

\noindent where $K$ is the number of MC samples, $\btheta^k$ is a random sample drawn from $p(\btheta)$, and $\mathbbm{1}_{\boldsymbol{\theta}^k \in \Theta (\mathbf{d})}$ is the indicator function. This function takes a value of one if  $\btheta^k$ is feasible (i.e., $\psi(\mathbf{d},\btheta^k)\leq 0$) or takes a value of zero otherwise. Note that $\mathbbm{1}_{\boldsymbol{\theta}^k \in \Theta (\mathbf{d})}=\mathbbm{1}_{\psi(\mathbf{d},\btheta^k)\leq 0}$.  From the law of large numbers we have that this approximation becomes exact as $|K| \rightarrow \infty$ \cite{shapiro2013sample}. 

For a given sample $\btheta^k$, problem \eqref{eq:feasibility_function_alt} can be written as:
\begin{equation}
\begin{aligned}
&\psi(\mathbf{d}, \boldsymbol{\theta}^k) = &\min_{\mathbf{z}^k, y^k} &&& y^kU \\
&&\text{s.t.} &&& f_j(\mathbf{d}, \mathbf{z}^k, \boldsymbol{\theta}^k) \leq y^kU, & j \in J \\
&&&&& h_i(\mathbf{d}, \mathbf{z}^k, \boldsymbol{\theta}^k) = 0, & i \in I\\
&&&&& y^k\in \{0,1\}
\end{aligned}
\label{eq:feasibility_function2}
\end{equation}
where $U\in \mathbb{R}_+$ is a sufficiently large constant. Note that the continuous variable for a particular sample $u^k$ is equivalently replaced with $y^kU$. The binary variable $y^k$ takes a value of zero if the realization $\btheta^k$ is feasible (corresponding to $\psi(\mathbf{d},\btheta^k)=0$) or takes a value of one if it is infeasible (corresponding to $\psi(\mathbf{d},\btheta^k)>0$).  

By combining all realizations $k\in K$, we obtain the aggregated problem:
\begin{equation}
\begin{aligned}
&&\min_{\mathbf{z}^k, y^k} &&& \frac{1}{|K|}\sum_{k \in K}y^kU \\
&&\text{s.t.} &&& f_j(\mathbf{d}, \mathbf{z}^k, \boldsymbol{\theta}^k) \leq y^k U, & j \in J,\; k\in K \\
&&&&& h_i(\mathbf{d}, \mathbf{z}^k, \boldsymbol{\theta}^k) = 0, & i \in I,\; k\in K \\
&&&&& y^k \in \{0,1\}, & k\in K.
\end{aligned}
\label{eq:inverse_SF}
\end{equation}
Since $\mathbf{d}$ is fixed, this problem is fully decoupled in the set $K$ (i.e., solving Problem \eqref{eq:inverse_SF} is equivalent to solving Problem \eqref{eq:feasibility_function2} for each $k\in K$) and the optimal values of the binary variables $y^k$ satisfy $SF_K(\mathbf{d})=\frac{1}{|K|}\sum_{k\in K}(1-y^k)$. Consequently, \eqref{eq:inverse_SF} can be used to compute the sample average approximation of the $SF$ index and the approximation becomes exact as $|K| \rightarrow \infty$. 

The constant $U$ in the objective can be eliminated without affecting the solution.  Moreover, we obtain an equivalent problem by maximizing $SF_K(\mathbf{d})=\frac{1}{|K|}\sum_{k \in K}(1-y^k)$ directly. These modifications give the problem:

\begin{equation}
\begin{aligned}
&&\max_{\mathbf{z}^k, y^k} &&& \frac{1}{|K|}\sum_{k \in K}(1-y^k) \\
&&\text{s.t.} &&& f_j(\mathbf{d}, \mathbf{z}^k, \boldsymbol{\theta}^k) \leq y^k U, & j \in J,\; k\in K \\
&&&&& h_i(\mathbf{d}, \mathbf{z}^k, \boldsymbol{\theta}^k) = 0, & i \in I,\; k\in K \\
&&&&& y^k \in \{0,1\}, & k\in K.
\end{aligned}
\label{eq:inverse_SF2}
\end{equation}

An important observation is that a continuous relaxation of problem \eqref{eq:inverse_SF}  is given by: 
\begin{equation}
\begin{aligned}
&&\min_{\mathbf{z}^k, u^k} &&& \frac{1}{|K|}\sum_{k \in K}u^k \\
&&\text{s.t.} &&& f_j(\mathbf{d}, \mathbf{z}^k, \boldsymbol{\theta}^k) \leq u^k, & j \in J,\; k\in K \\
&&&&& h_i(\mathbf{d}, \mathbf{z}^k, \boldsymbol{\theta}^k) = 0, & i \in I,\; k\in K \\
&&&&& u^k \geq 0, & k\in K.
\end{aligned}
\label{eq:inverse_SF3}
\end{equation}
This provides a relaxation because $u^k$ can be modeled as $u^k=y^kU$ with $0\leq y^k\leq 1$ and sufficiently large $U\in \mathbb{R}_+$. In summary, the continuous problem \eqref{eq:inverse_SF3} is a {\em relaxation} of the mixed-integer problem \eqref{eq:inverse_SF} (or, equivalently, \eqref{eq:inverse_SF2}). We also note that $u_k=\psi(\mathbf{d},\btheta^k)$ is a measure of infeasibility for sample $\btheta^k$ and thus the continuous problem \eqref{eq:inverse_SF3} minimizes the {\em mean (expected) infeasibility}.  In other words, the expected infeasibility is a surrogate measure of flexibility. 
 
\subsection{Mixed-Integer Optimal Design Problem}

We now proceed to find a design $\mathbf{d}$ that minimizes the design cost and that maximizes the stochastic flexibility index.  A sample average approximation of this problem  is given by:
\begin{equation}
\begin{aligned}
& &\min_{\mathbf{d}\in \mathcal{D}, \mathbf{z}^k, y^k} &&& \left\{c(\mathbf{d}), -\frac{1}{|K|}\sum_{k \in K}(1-y^k)\right\} \\
&&\text{s.t.} &&& f_j(\mathbf{d}, \mathbf{z}^k, \boldsymbol{\theta}^k) \leq y^k U, && j \in J, \ \ k \in K \\
&&&&& h_i(\mathbf{d}, \mathbf{z}^k, \boldsymbol{\theta}^k) = 0, && i \in I, \ \ k \in K \\
&&&&& y^k\in \{0,1\}, && k \in K \\
\end{aligned}
\label{eq:2_stage}
\end{equation}
A Pareto solution of this problem can be obtained by solving an $\epsilon$-constrained problem of the form:
\begin{equation}
\begin{aligned}
& &\min_{\mathbf{d}\in \mathcal{D}, \mathbf{z}^k, y^k} &&&  c(\mathbf{d})\\
&&\text{s.t.} &&& f_j(\mathbf{d}, \mathbf{z}^k, \boldsymbol{\theta}^k) \leq y^k U, && j \in J, \ \ k \in K \\
&&&&& h_i(\mathbf{d}, \mathbf{z}^k, \boldsymbol{\theta}^k) = 0, && i \in I, \ \ k \in K \\
&&&&& \frac{1}{|K|}\sum_{k \in K}(1-y^k)\geq \epsilon_s\\
&&&&& y^k\in \{0,1\}&& k \in K \\
\end{aligned}
\label{eq:2_stageeps}
\end{equation}
where $\epsilon_s$ is the $\epsilon$-parameter. This formulation is precisely a sample average approximation of the joint chance-constrained problem \eqref{eq:design}. This problem is computationally expensive to solve due to its mixed-integer nature and due to the large number of potential MC samples. 

\subsection{Continuous Optimal Design Problem} \label{sec:cont_formulation}
A key observation is that Pareto solutions for the conflict resolution problem can also be obtained by solving the mixed-integer problem:
\begin{equation}
\begin{aligned}
& &\max_{\mathbf{d}\in \mathcal{D}, \mathbf{z}^k, y^k} &&& \frac{1}{|K|}\sum_{k \in K}(1-y^k) \\
&&\text{s.t.} &&& f_j(\mathbf{d}, \mathbf{z}^k, \boldsymbol{\theta}^k) \leq y^k U, && j \in J, \ \ k \in K \\
&&&&& h_i(\mathbf{d}, \mathbf{z}^k, \boldsymbol{\theta}^k) = 0, && i \in I, \ \ k \in K \\
&&&&&c(\mathbf{d})\leq \epsilon_c\\
&&&&& y^k\in \{0,1\}, && k \in K.
\end{aligned}
\label{eq:2_stageeps_invert}
\end{equation}
Note that this is simply the counterpart to \eqref{eq:2_stageeps}. We denote the solution of this problem as $\mathbf{d}^*,{y^k}^*$ and we have that $c(\mathbf{d}^*)= \epsilon_f$  holds at the solution (because the objective is conflicting) and the Pareto pair is given by $(c(\mathbf{d}^*),SF_K)=(\epsilon_f,SF_K)$, where $ SF_K=\frac{1}{|K|}\sum_{k \in K}(1-{y^k}^*)$. 

A continuous relaxation of this problem is given by:
\begin{equation}
\begin{aligned}
& &\max_{\mathbf{d}\in \mathcal{D}, \mathbf{z}^k, y^k} &&& \frac{1}{|K|}\sum_{k \in K}(1-y^k) \\
&&\text{s.t.} &&& f_j(\mathbf{d}, \mathbf{z}^k, \boldsymbol{\theta}^k) \leq y^kU, && j \in J, \ \ k \in K \\
&&&&& h_i(\mathbf{d}, \mathbf{z}^k, \boldsymbol{\theta}^k) = 0, && i \in I, \ \ k \in K \\
&&&&&c(\mathbf{d})\leq \epsilon_c\\
&&&&& 0\leq y^k\leq 1, && k \in K.
\end{aligned}
\label{eq:2_stageeps2}
\end{equation}
The solution of this problem delivers a continuous solution $\bar{y}^k$ from which we recover integer values as $\bar{y}^k\leftarrow \mathbbm{1}_{\bar{y}^k}$ (i.e., we recover $\bar{y}^k = 0$ if $k$ is feasible or $\bar{y}^k = 1$ otherwise). With this, we compute an estimate of the SF index as $\bar{SF}_K=\frac{1}{|K|}\sum_{k \in K}(1-\mathbbm{1}_{\bar{y}^k})$. We solve \eqref{eq:2_stageeps2} again (by adjusting $\mathbf{d}\in \mathcal{D}$ and $\mathbf{z}^k$) by fixing the binary variables to ensure that the selection of {\em active} constraints provides a feasible solution. Note that the objective function $\bar{SF}_K$ is fixed since the binaries are fixed. The solution of this feasibility problem is given by $\bar{\mathbf{d}},\bar{\mathbf{z}}^k$. We highlight that problem \eqref{eq:2_stageeps2} can also be solved by using the continuous variables $u^k$. We will see that this approach delivers Pareto pairs $(c(\bar{\mathbf{d}}),\bar{SF}_K)$ that closely match the optimal Pareto pairs $(c({\mathbf{d}^*}),{SF}_K)$ of the mixed-integer formulation. We also note an approximation of comparable quality is achieved when the counterpart problem of \eqref{eq:2_stageeps2} is solved.

A precise theoretical justification for this continuous relaxation providing such a high quality approximation is the subject of ongoing research and is beyond the scope of this work. We hypothesize that the high quality of the solutions is due to the degenerate nature of the SF index. Specifically, different combinations of active constraints give the same or a very similar SF index. This behavior has been recently reported in \cite{curtis2018sequential}, where the authors note that chance constraints give rise to a wide range of local minima with similar values. This behavior becomes more evident when one moves the SF index into the objective (as opposed to imposing it in the constraints). This is precisely why it is important to think about the design problem as a conflict resolution problem. Recent work has also provided evidence that continuous (Lagrangian) relaxations of chance-constrained problems provide high quality approximations \cite{ahmed2017nonanticipative}.

\section{Case Studies and Analysis} \label{sec:cases}

We analyze the behavior of the proposed formulations by applying them to distribution networks. We consider a simple three-node network, the IEEE-14 power distribution network, and a large 141-node network. All formulations are implemented in JuMP 0.18.5 \cite{DunningHuchetteLubin2017} and are solved using CPLEX 12.6.3 for mixed-integer problems and Gurobi 7.5.1 for linear problems on a dual Intel\textregistered \, Xeon\textregistered \, ES-2698 v3 machine running at 2.30 GHz with 64 hardware threads and 198 GB of RAM. All results can be reproduced using the scripts provided in \url{https://github.com/zavalab/JuliaBox/tree/master/FlexDesign}. 

\subsection{Case Study Models}
Each network is modeled by performing balances at each node $n \in \mathcal{C}$ and enforcing capacity constraints on the arcs $a_l, l \in \mathcal{A},$ and on the suppliers $s_b, b \in \mathcal{S}$. The demands $r_m, m \in \mathcal{R},$ are assumed to be the uncertain parameters which are given by the MC samples. The deterministic network model is given by:
\begin{subequations}
	\begin{equation}
	\sum_{l \in \mathcal{A}_n^{rec}} a_l - \sum_{l \in \mathcal{A}_n^{snd}} a_l + \sum_{b \in \mathcal{S}_n} s_b - \sum_{m \in \mathcal{R}_n} r_m = 0, \ \ \ n \in \mathcal{C}
	\label{eq:node_balances}
	\end{equation}
	\begin{equation}
	-a_l^C - d^a_l \leq a_l \leq a_l^C + d^a_l, \ \ \ l \in \mathcal{A}
	\label{eq:arc_constrs}
	\end{equation}
	\begin{equation}
	0 \leq s_b \leq s_b^C + d^s_b, \ \ \ b \in \mathcal{S}
	\label{eq:supp_constrs}
	\end{equation}
	\label{eq:network_model}
\end{subequations}
\noindent where $\mathcal{A}_n^{rec}$ denotes the set of receiving arcs at node $n$, $\mathcal{A}_n^{snd}$ denotes the set of sending arcs at $n$, $\mathcal{S}_n$ denotes the set of suppliers at $n$, $\mathcal{R}_n$ denotes the set of demands at $n$, $a_l^C$ are the arc capacities, $d^a_l \geq 0$ are design variables that increase arc capacity, $s_b^C$ are the supplier capacities, and $d^s_b \geq 0$ are design variables that increase supplier capacity. 

\begin{figure}[!htb]
	\centering
	\includegraphics[width=0.4\textwidth]{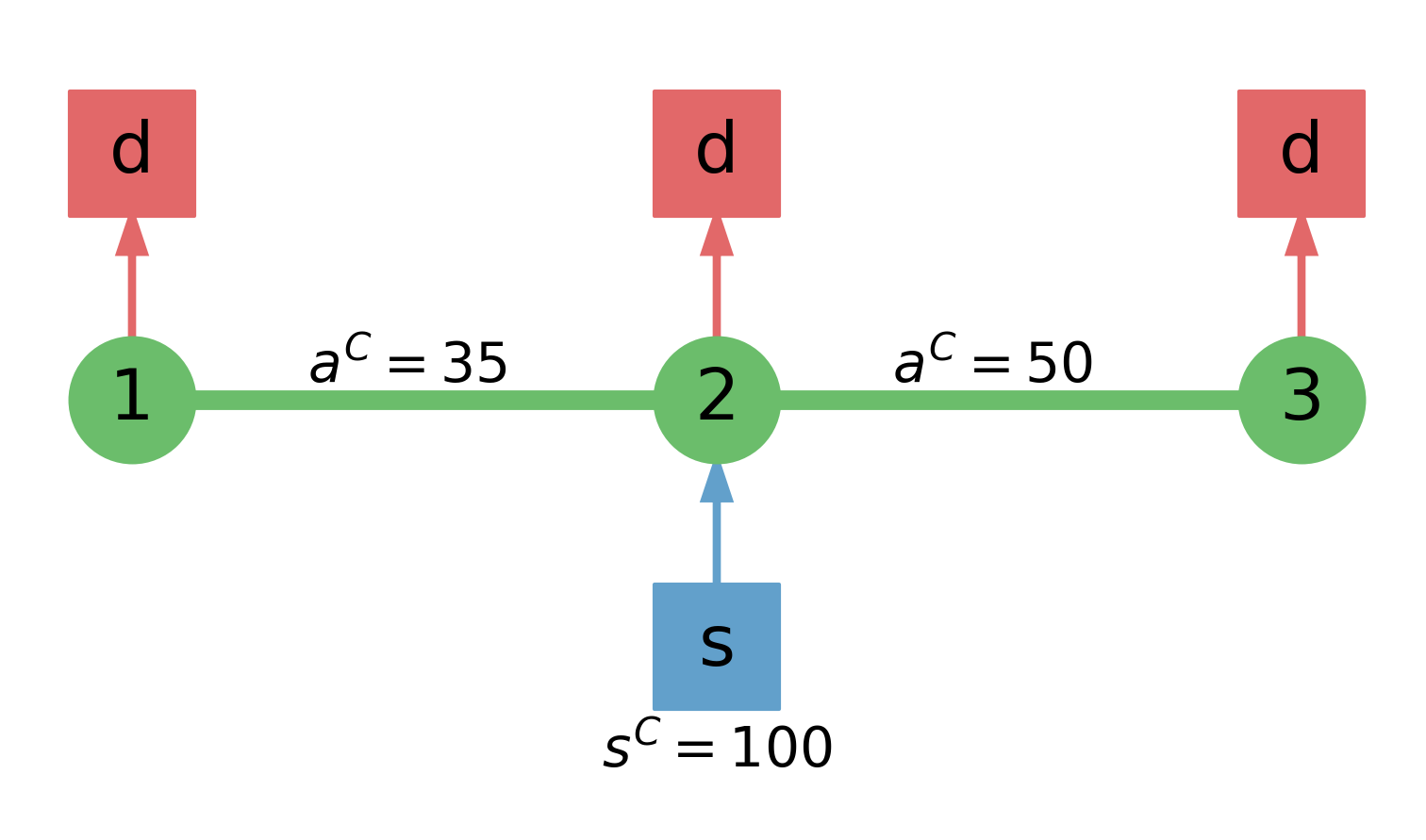}
	\caption{The three-node distribution network}
	\label{fig:3d_design}
\end{figure}

The three-node distribution network features a centralized supplier configuration. Figure \ref{fig:3d_design} details this network and provides the arc and supplier capacities. The network is subjected to multivariate Gaussian demands $\btheta = (r_1,r_2,r_3) \sim \mathcal{N}(\bar{\btheta}, V_{\btheta})$. The mean $\bar{\btheta}$ is taken to be $\bar{\btheta}=(0.0, 60.0, 10.0)$, and the covariance matrix $V_{\btheta}$ is assumed to be:

\begin{equation}
V_{\btheta} = 
\begin{bmatrix}
80 & 0 & 0 \\
0 & 80 & 0 \\
0 & 0 & 120
\end{bmatrix}.
\label{eq:3d_covars}
\end{equation}
We also set the parameter $U = 10000$.

\begin{figure}[!htb]
	\centering
	\includegraphics[width=0.4\textwidth]{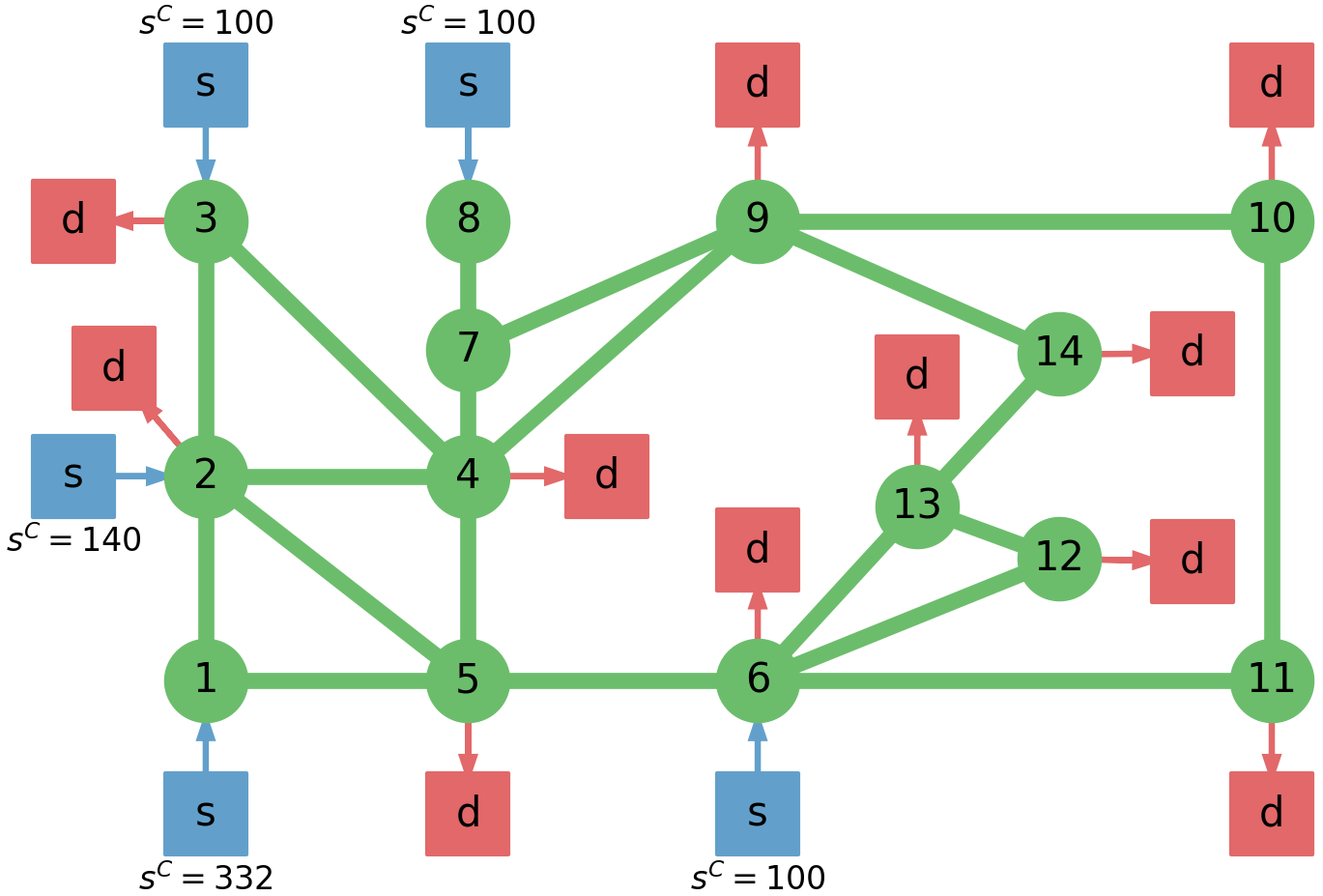}
	\caption{Schematic of the IEEE 14-node power system where the values $s^C$ are indicated and all the values $a^C = 100$.}
	\label{fig:14_diagram}
\end{figure}

The IEEE 14-node network was originally provided in Dabbagchi in \cite{ieee14_origin}. The system data is obtained from MATPOWER. This test case does not provide arc capacities so we enforce $a^C = 100$ for all the arcs. A schematic of this system is provided in Figure \ref{fig:14_diagram}. This system is subjected to a total of 11 uncertain parameters (the network demands). The demands are assumed to be $\btheta \sim \mathcal{N}(\bar{\btheta}, V_{\btheta})$, where $\bar{\btheta} = (87.3,50.0,25.0,28.8,50.0,25.0,\allowbreak0,0,0,0,0)$ and $V_{\btheta}$ is symmetric matrix with $(V_{\btheta})_{ii} = 1200$ and $(V_{\btheta})_{ij} = 240, \ \forall i \neq j$. Also, we set $U = 10000$.

\begin{figure}[!htb]
	\centering
	\includegraphics[width=0.8\textwidth]{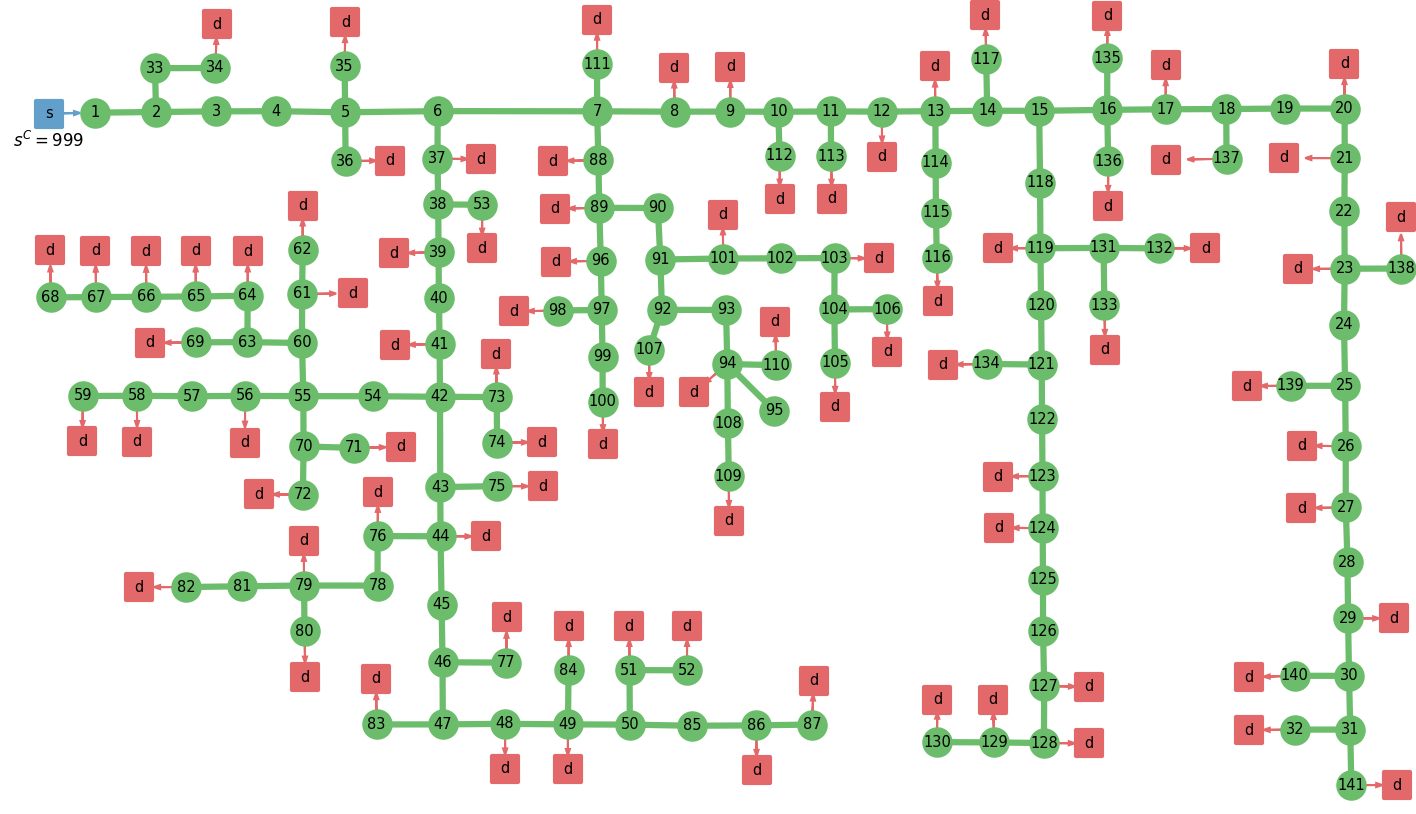}
	\caption{Schematic of the 141-node power distribution network.}
	\label{fig:141_diagram}
\end{figure}

The 141-node power distribution network corresponds to an urban area in Caracas, Venezuela and was originally developed in \cite{khodr2008maximum}. The network data is extracted from MATPOWER, but again no arc capacities are provided (we assume $a^C = 100$). Figure \ref{fig:141_diagram} provides a schematic of the network. This system is subjected to 84 uncertain disturbances corresponding to the demands. The demands are assumed to be $\btheta \sim \mathcal{N}(\bar{\btheta}, V_{\btheta})$, where $\bar{\btheta} = \bar{\btheta}_{fc}$ and $V_{\btheta} = 100 \mathbb{I}$. The point $\btheta_{fc}$ is the feasible center and is the instance of $\btheta \in \Theta(\mathbf{d})$ that minimizes the feasibility function $\psi(\mathbf{d}, \btheta)$ \cite{pulsipher2018rank}. We set the upper bounding constant to $U = 10000$. 

For each network, the cost function $c(\mathbf{d})$ is defined to be linear of the form:
	 \begin{equation}
	 c(\mathbf{d}^s, \mathbf{d}^a) = \frac{1}{\sqrt{n_d}} \left(\sum_{b \in \mathcal{S}} d^s_b +  \sum_{l \in \mathcal{A}} d^a_l \right).
	 \end{equation}
Here the unit cost for each design variable is taken to be $1/\sqrt{n_d}$, where $n_d = |\mathcal{S}| + |\mathcal{A}|$. Furthermore, a relatively large value of $U$ is selected in the above problems to ensure a sufficient amount of slack is provided for the inequality constraints. This is all done for convenience in conducting the analysis below.

\FloatBarrier 

\subsection{Mixed-Integer Formulation} \label{sec:bin_2stage}
Combining \eqref{eq:2_stageeps_invert} with the model equations in \eqref{eq:network_model}, we obtain:

\begin{equation}
\begin{aligned}
& &\max_{d^s_b, d^a_l, a_l^k, s_b^k, y^k} &&& \frac{1}{|K|} \sum_{k \in K} (1-y^k) \\
&&\text{s.t.} &&& -a_l^C - d^a_l - a_l^k \leq y^k U, && l \in \mathcal{A}, \ \ k \in K \\
&&&&& -a_l^C - d^a_l + a_l^k \leq y^k U, && l \in \mathcal{A}, \ \ k \in K \\
&&&&& -s_b^k \leq y^k U, && b \in \mathcal{S}, \ \ k \in K \\
&&&&& -s_b^C - d^s_b + s_b^k \leq y^k U, && b \in \mathcal{S}, \ \ k \in K \\
&&&&& \sum_{l \in \mathcal{A}_n^{rec}} a_l^k - \sum_{l \in \mathcal{A}_n^{snd}} a_l^k + \sum_{b \in \mathcal{S}_n} s_b^k - \sum_{m \in \mathcal{R}_n} r_m^k = 0, && n \in \mathcal{C}, \ \ k \in K \\
&&&&& \frac{1}{\sqrt{n_d}} \left(\sum_{b \in \mathcal{S}} d^s_b +  \sum_{l \in \mathcal{A}} d^a_l \right) \leq \epsilon_c \\
&&&&& y^k \in \{0,1\} && k \in K \\
&&&&& d^s_b \geq 0, \ \ d^a_l \geq 0,&& b \in \mathcal{S}, \ \  l \in \mathcal{A}.
\end{aligned}
\label{eq:2_stage_dist}
\end{equation}

\begin{figure}[!htb]
	\centering
		\includegraphics[width=0.6\textwidth]{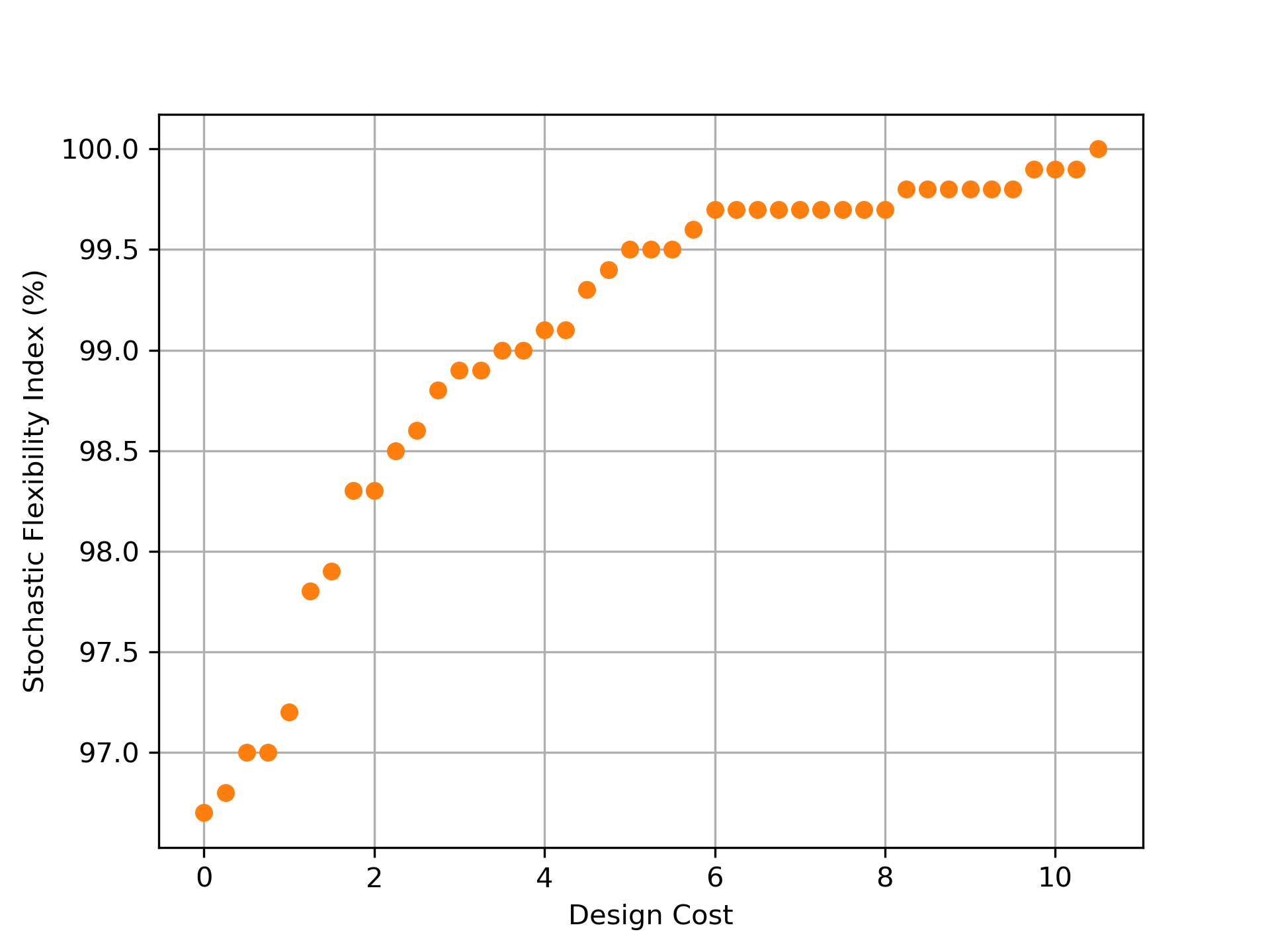}
	\caption{The Pareto set for the three-node network obtained with the mixed-integer optimal design formulation.}
	\label{fig:3d_pareto_bin}
\end{figure}

\noindent We extract elements of the Pareto set by varying the cost threshold $\epsilon_c$. This analysis is done using the three-node network with 1,000 MC samples that are generated from the underlying distribution. The Pareto set is obtained from Problem \eqref{eq:2_stage_dist} for 43 values of $\epsilon_c$ set from 0 to 10.5 in 0.25 increments. A subset of the numerical results is presented in Table \ref{tab:3d_bin_sub} and the complete set of results is given in Tables \ref{tab:3d_data1} and \ref{tab:3d_data2} in the Appendix. The minimum possible $SF_K$ index associated with the base case is 96.7\%. Initially, as we increase $\epsilon_c$, we are able to increase the $SF_K$ index at a relatively small cost. For instance, a 2.1\% improvement of the $SF_K$ index only incurs a design cost of 2.75 (1.31 per \% increase). In contrast, the final 0.4\% increase in the $SF_K$ index, making the network perfectly flexible relative to the 1,000 sampled instances, incurs a 4.75 increase in design cost (11.88 per \% increase). This trend can be visualized in Figure \ref{fig:3d_pareto_bin}.

\begin{table}[!htb]
	\caption{A reduced subset of results obtained for the 3-node network with Problem \eqref{eq:2_stage_dist} using 1,000 MC samples}.
	\begin{center}
		\begin{tabular}{|c|ccc|}
			\hline
			$\epsilon_c$ & Design Cost & $SF_K$ (\%) & Solution Time (s) \\ \hline \hline
			0    & 0     & 96.7  & 0.0155 \\
			1.25 & 1.25  & 97.8  & 0.0190 \\
			2.75 & 2.75  & 98.8  & 0.0175 \\
			4.25 & 4.25  & 99.1  & 0.0200 \\
			5.75 & 5.75  & 99.6  & 0.0193 \\
			7.25 & 7.25  & 99.7  & 0.0207 \\
			8.75 & 8.75  & 99.8  & 0.0192 \\
			10.5 & 10.5  & 100   & 0.0195 \\ \hline
		\end{tabular}
	\end{center}
	\label{tab:3d_bin_sub}
\end{table}

\begin{figure}[!htb]
	\centering
	\begin{subfigure}[t]{.31\textwidth}
		\centering
		\includegraphics[width=\textwidth]{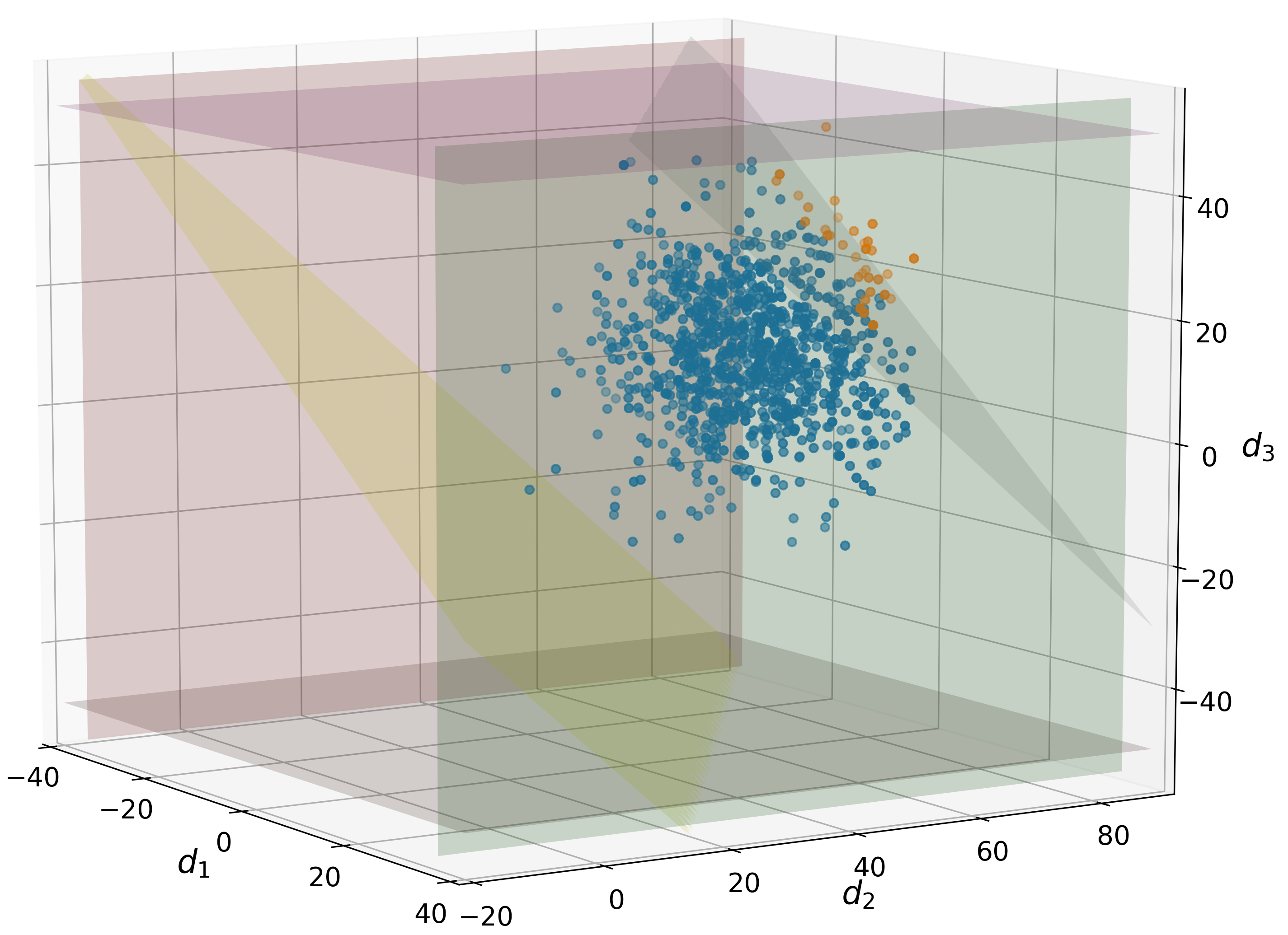}
		\caption{$\epsilon_c = 0$}
	\end{subfigure}
	\quad
	\begin{subfigure}[t]{.31\textwidth}
		\centering
		\includegraphics[width=\textwidth]{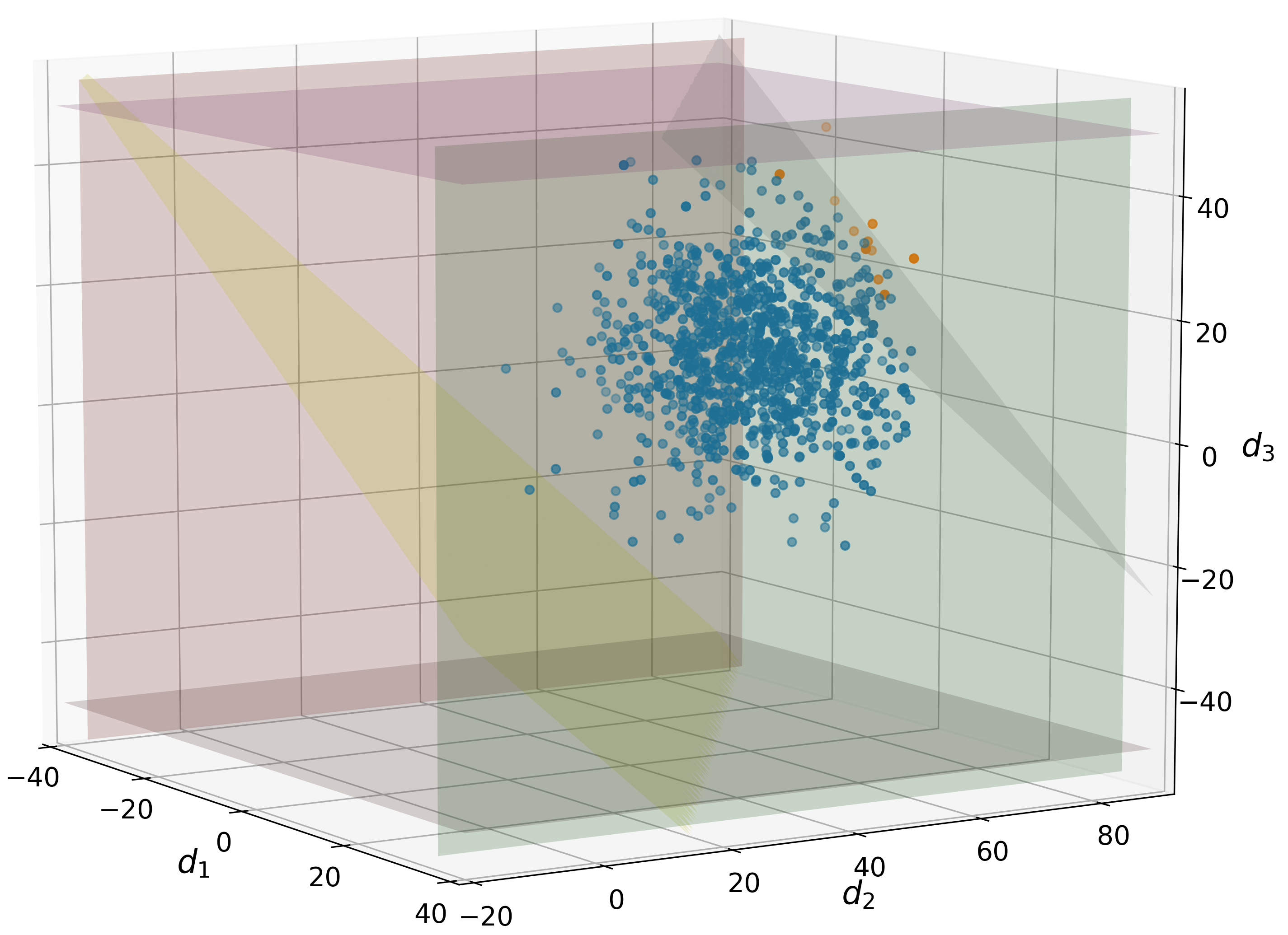}
		\caption{$\epsilon_c = 2.75$}
	\end{subfigure}
	\quad
	\begin{subfigure}[t]{.31\textwidth}
		\centering
		\includegraphics[width=\textwidth]{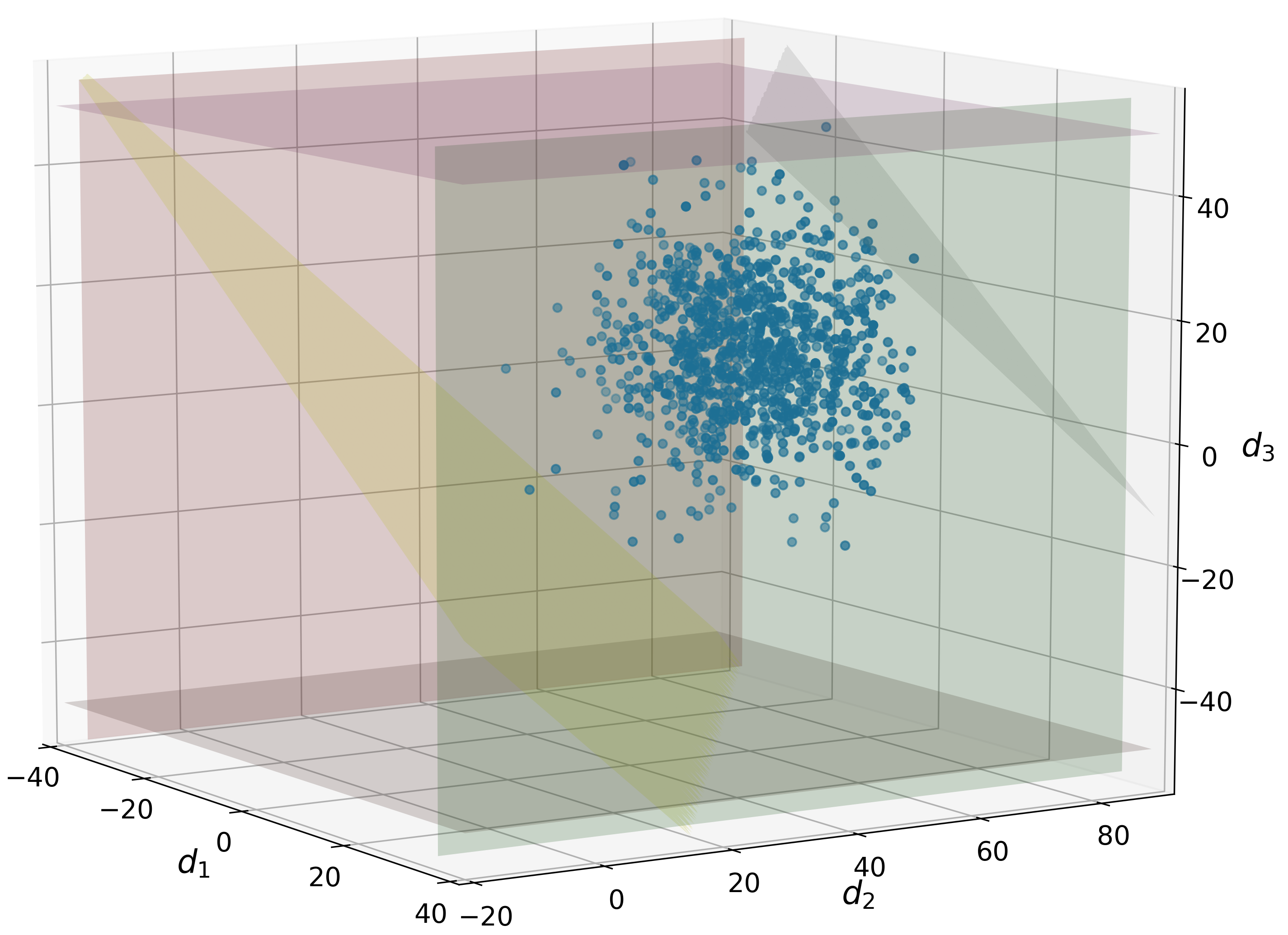}
		\caption{$\epsilon_c = 10.5$}
	\end{subfigure}
	\caption{The optimized solutions to Problem \eqref{eq:2_stage_dist} for the three-node network, showing how the design constraints are shifted to minimize the number of infeasible (orange) instances as the value of $\epsilon_c$ is increased.}
	\label{fig:3d_plotted_solns}
\end{figure}

We substitute the node balances into the capacity constraints to obtain a system of inequalities that are expressed solely in terms of $\btheta$. Thus, the three-node network solutions can be visualized in three dimensions. Figure \ref{fig:3d_plotted_solns} shows three solutions corresponding to $\epsilon_c = \{0, 2.75, 10.5\}$. Here, we observe how the diagonal plane (on the right side of the figures), which corresponds the supplier capacity constraint, is shifted as the value of $\epsilon_c$ is increased to minimize the number of infeasible instances and thus maximize the $SF_K$ index. 

The same analysis is done using the IEEE-14 network with 2,000 MC samples. This creates a MILP with 50,025 continuous variables, 2,000 binary variables, and 128,001 constraints. The Pareto set is created by varying the value of $\epsilon_c$ from 0 to 65 in increments of 2.5. All of the corresponding numerical results are given in Table \ref{tab:ieee14_data} in the appendix. Figure \ref{fig:ieee14_pareto_bin} shows the Pareto set corresponding to these results. Here, the minimum $SF_K$ index associated with the base case is 88.10\% and the maximum $SF_K$ index is 91.50\%. This means that, regardless of how much the line and supplier capacities are increased, we can only improve $SF_K$ by 3.4\% (relative to this MC sample set). This behavior occurs because in certain sampled instances a number of demands are sufficiently negative such that this is a surplus in the system that cannot be rectified since there are no sinks in the network. Such behavior is not readily obvious and thus highlights the utility of this optimal design framework in understanding how to promote the flexibility of a system and in determining to what extent its flexibility can be improved. We also observe that the average solution time for the MILP problems  is 40.47s.

\begin{figure}[!htb]
	\centering
	\includegraphics[width=0.6\textwidth]{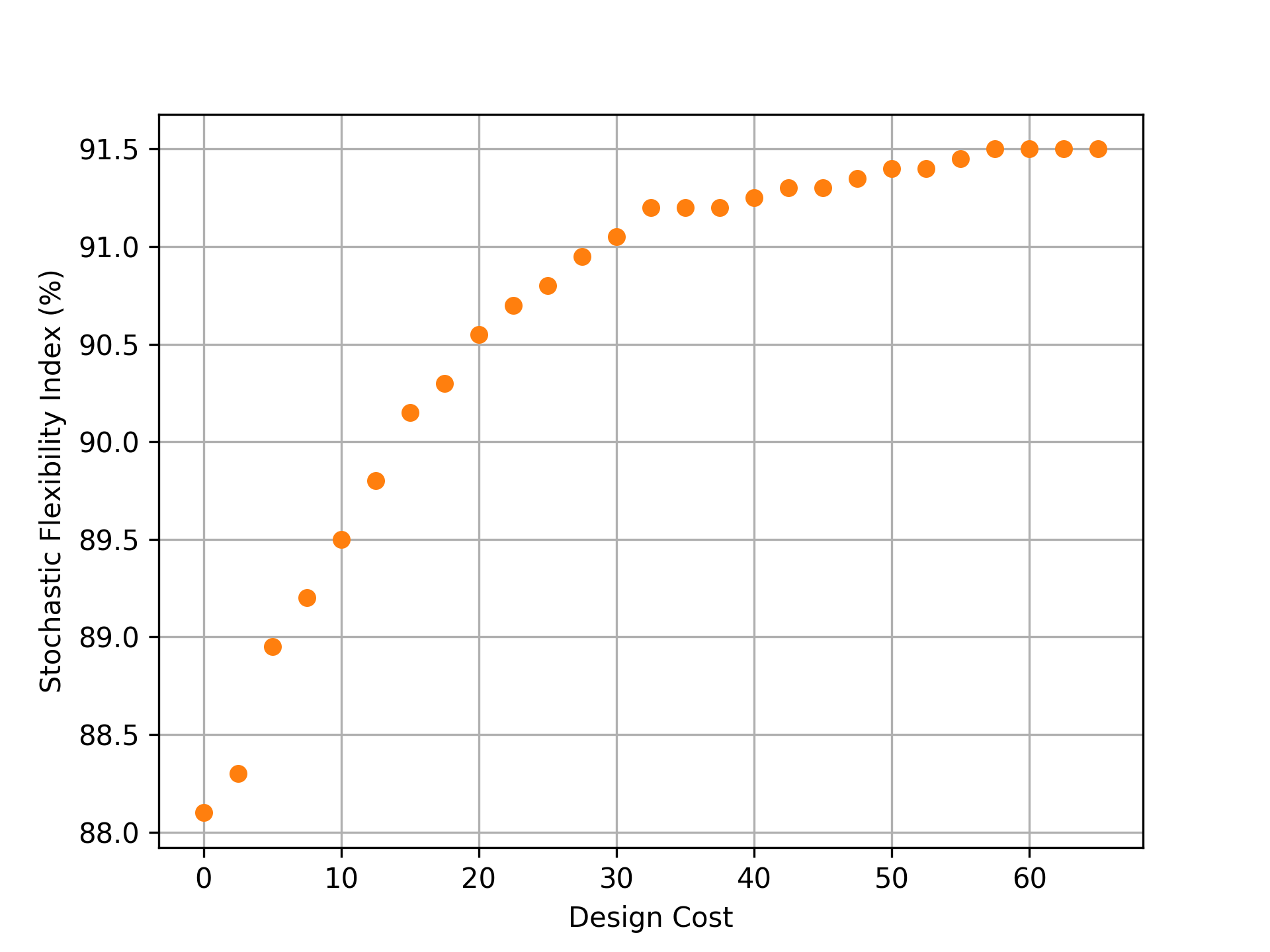}
	\caption{The Pareto set for the IEEE-14 power network obtained with the mixed-integer design formulation.}
	\label{fig:ieee14_pareto_bin}
\end{figure}

Problem \eqref{eq:2_stage_dist} is also applied to the 141-node power distribution network using 10,000 MC samples and is evaluated for each value of $\epsilon_c$ from 0 to 300 in increments of 10. Each MILP contains 1,410,141 continuous variables, 10,000 binary variables, and 4,230,001 constraints. A time limit of 3,600s is also imposed. All 30 results are summarized in Table \ref{tab:141_data} in the Appendix and Figure \ref{fig:141_pareto_bin} shows these results. Here, a much more appreciable range of $SF_K$ index values is obtained where the minimum $SF$ index of the system is 35.60\% and the maximum $SF$ index is 70.64\%. This limited maximum $SF$ index can be attributed to topological limitations in the network. Specifically, the sampled demands are Gaussian random variables and thus can be negative which in certain cases means the network is confronted with excess power production that it cannot shed. Thus, we observe that in this case increasing the network capacities alone cannot provide a $SF$ index near 100\% due to topographical limitations. Instead, the system topology would need to be modified to enhance its flexibility (e.g., add a sink to shed excess power production). This highlights how this framework is useful in guiding system design by demonstrating to what degree the system flexibility can be bettered via continuous design variables such as capacity.

A key observation is that only 12 of the 30 problems converged within the time limit. This unfavorable behavior is likely due to the large-U constraints which can lead to weak linear program (LP) relaxations. This scalability limitation clearly demonstrates that Problem \eqref{eq:2_stageeps_invert} can become impractical for moderate and large sized systems. We note some advanced methods such as those proposed in \cite{song2014chance,qiu2014covering} can be used to select smaller $U$ values that might strengthen the LP relaxations. However, our proposed continuous relaxation entails a more efficient strategy since it does not require the solution of a MIP. Another observation drawn from the solution times in Table \ref{tab:141_data} is that Pareto pairs near minimum or maximum design cost require significantly less computational time relative to the intermediate pairs. This trend is also prevalent in the computational results of the IEEE-14 network and the three-node network.

\begin{figure}[!htb]
	\centering
	\includegraphics[width=0.6\textwidth]{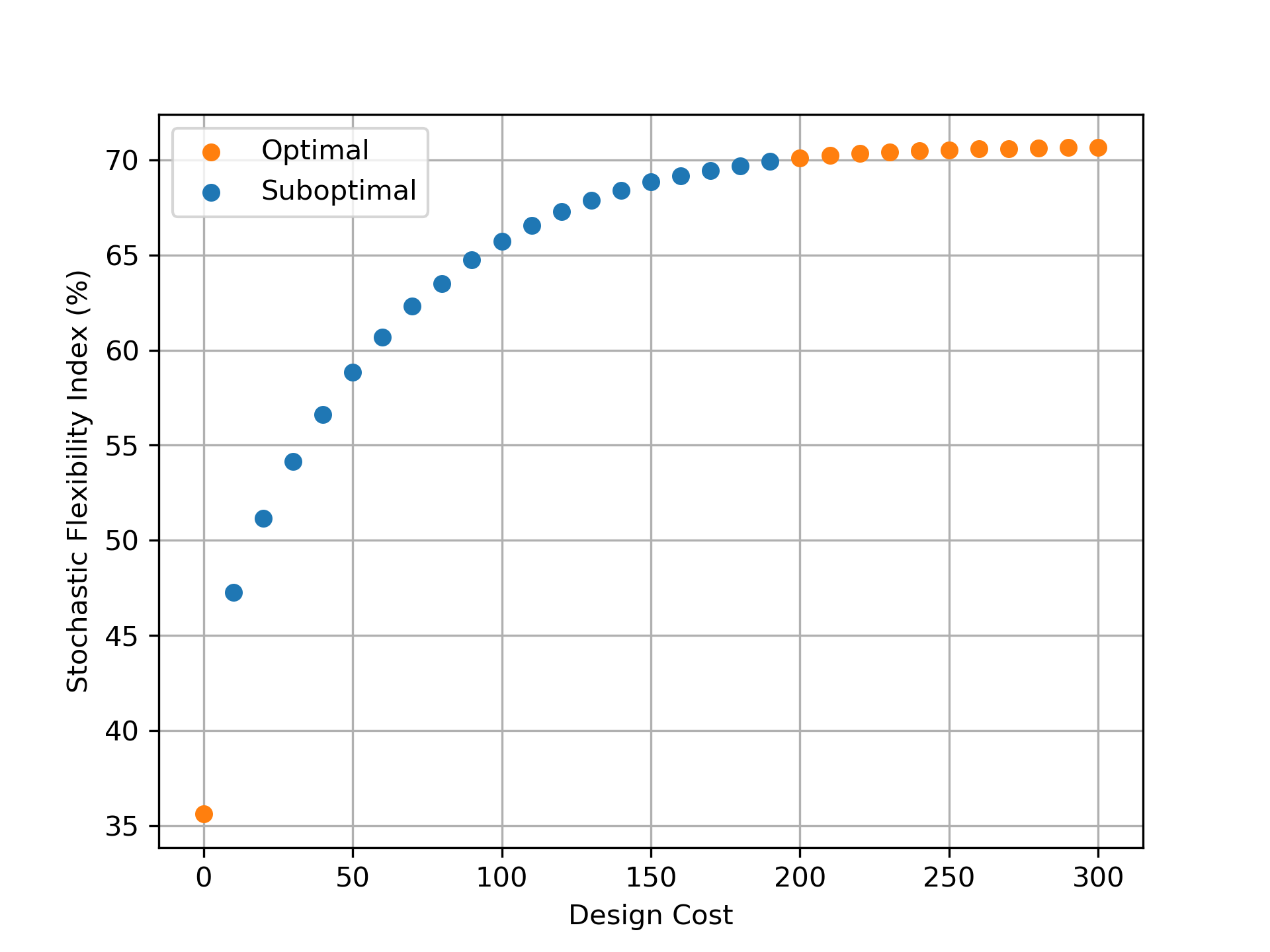}
	\caption{The Pareto set for the 141-node power network obtained with the mixed-integer design formulation. Here, the optimal points are those that solved within the time limit of 3,600s.}
	\label{fig:141_pareto_bin}
\end{figure} 

\FloatBarrier

\subsection{Continuous Formulation} \label{sec:cont_2stage}
We compare the Pareto pairs of the mixed-integer formulation to those obtained with the continuous formulation:
\begin{equation}
\begin{aligned}
& &\max_{d^s_b, d^a_l, a_l^k, s_b^k, y^k} &&& \frac{1}{|K|} \sum_{k \in K} (1-y^k) \\
&&\text{s.t.} &&& -a_l^C - d^a_l - a_l^k \leq y^k U, && l \in \mathcal{A}, \ \ k \in K \\
&&&&& -a_l^C - d^a_l + a_l^k \leq y^k U, && l \in \mathcal{A}, \ \ k \in K \\
&&&&& -s_b^k \leq y^k U, && b \in \mathcal{S}, \ \ k \in K \\
&&&&& -s_b^C - d^s_b + s_b^k \leq y^k U, && b \in \mathcal{S}, \ \ k \in K \\
&&&&& \sum_{l \in \mathcal{A}_n^{rec}} a_l^k - \sum_{l \in \mathcal{A}_n^{snd}} a_l^k + \sum_{b \in \mathcal{S}_n} s_b^k - \sum_{m \in \mathcal{R}_n} r_m^k = 0, && n \in \mathcal{C}, \ \ k \in K \\
&&&&& \frac{1}{\sqrt{n_d}} \left(\sum_{b \in \mathcal{S}} d^s_b +  \sum_{l \in \mathcal{A}} d^a_l \right) \leq \epsilon_c \\
&&&&& 0 \leq y^k \leq 1, && k \in K \\
&&&&& d^s_b \geq 0, \ \ d^a_l \geq 0,&& b \in \mathcal{S}, \ \  l \in \mathcal{A}.
\end{aligned}
\label{eq:2_stage_mod_dist}
\end{equation}

The $\bar{SF}_K$ index is determined after solving this problem via the indicator function $\mathbbm{1}_{y^k}$ as described in Section \ref{sec:cont_formulation}. The three-node network is again analyzed under the same conditions described in Section \ref{sec:bin_2stage} (i.e., the same samples and $\epsilon_c$ values). For the small 3-node network, the problem is an LP with 4,003 variables and 9,001 constraints. Each solution is also verified to be integer feasible by ensuring the formulation is feasible with the values of the $y^k$ variables fixed to their indicated values as described in Section \ref{sec:cont_formulation}. The numerical results are provided in Tables \ref{tab:3d_data1} and \ref{tab:3d_data2} in the Appendix and Figure \ref{fig:3d_pareto_cont} juxtaposes these results with those obtained with the mixed-integer formulation. Here, we note that the Pareto pairs are equivalent and that the mixed-integer $y^k$ values exactly match the indicated $y^k$ values obtained with the continuous formulation (i.e., the same active constraints are identified). Also, each continuous formulation solution only required 0.0049s on average to converge which equates to a 74\% reduction in computational cost. 

\begin{figure}[!htb]
	\centering
	\includegraphics[width=0.6\textwidth]{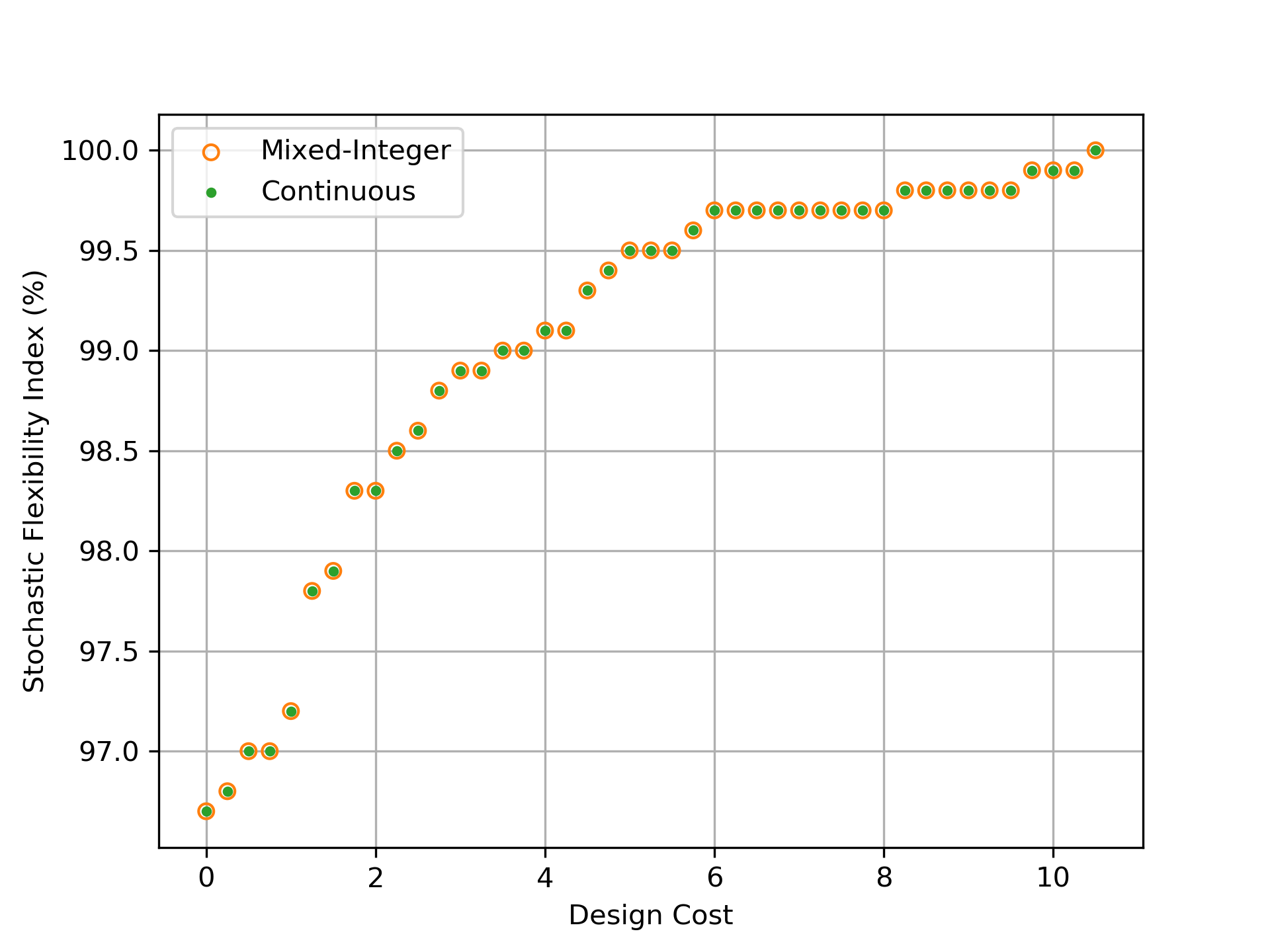}
	\caption{The Pareto sets for three-node network obtained with mixed-integer and continuous formulations.}
	\label{fig:3d_pareto_cont}
\end{figure}

We also applied Problem \eqref{eq:2_stage_mod_dist} to the IEEE-14 node power network with the same samples and conditions used in Section \ref{sec:bin_2stage}. This is an LP with 52,025 variables and 128,001 constraints. The indicated values of $y^k$ are all verified to be integer feasible. These results are detailed in Table \ref{tab:ieee14_data} in the Appendix and Figure \ref{fig:ieee14_pareto_cont} compares the results with those obtained with the mixed-integer formulation. Again, we note that the continuous formulation is able to recover the same Pareto set as the mixed-integer counterpart. Only the last Pareto pair differs (by one $y^k$ out of 2,000). All other pairs are identical and identify exactly the same sets of active constraints. Furthermore, the continuous formulation only requires 1.74s on average to solve, which corresponds to a 96\% reduction in computational cost relative to the mixed-integer formulation.

\begin{figure}[!htb]
	\centering
		\includegraphics[width=0.6\textwidth]{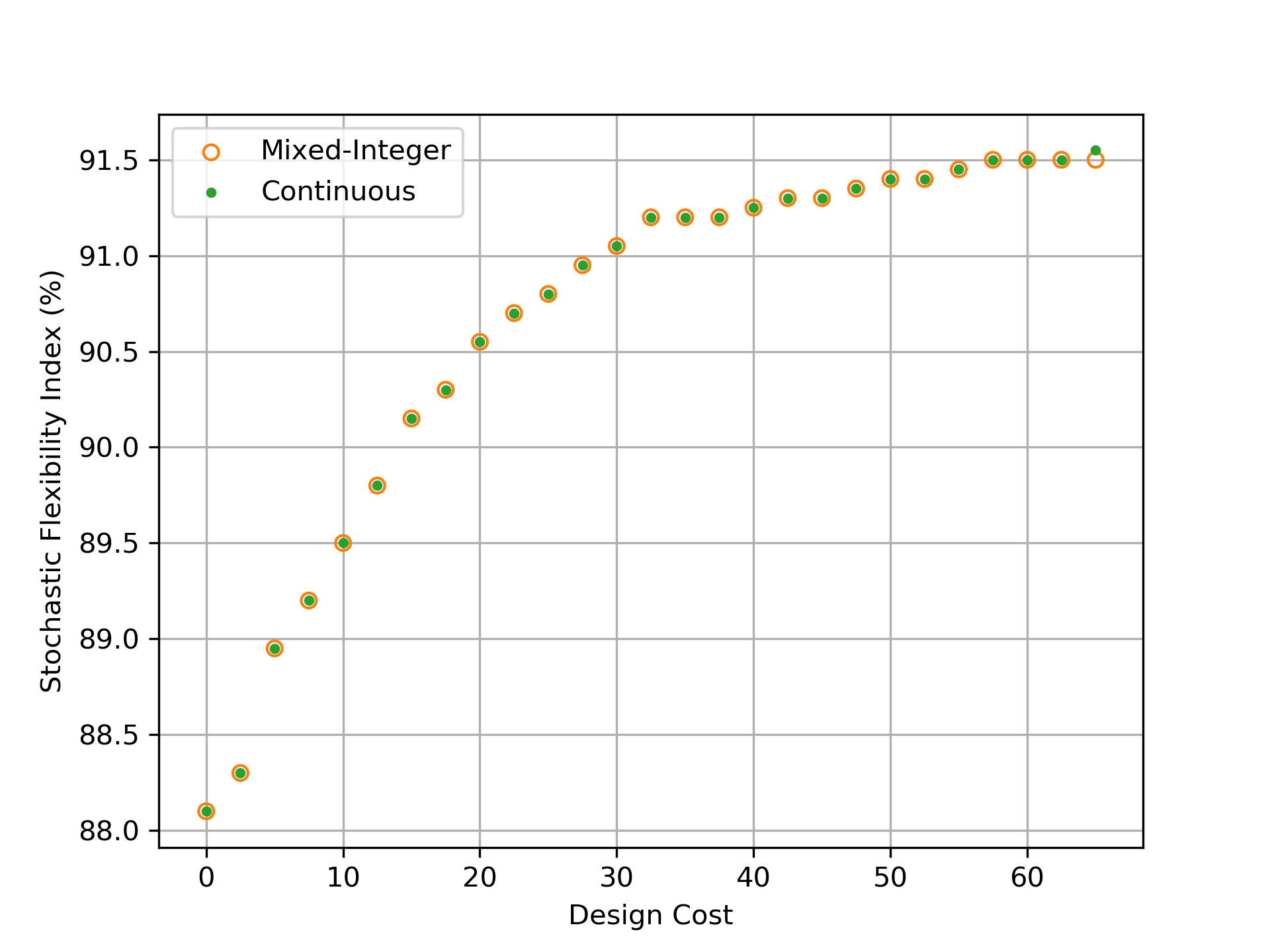}
	\caption{Pareto sets for IEEE-14 node power network obtained with mixed-integer and continuous formulations.}
	\label{fig:ieee14_pareto_cont}
\end{figure}

Finally, we applied Problem \eqref{eq:2_stage_mod_dist} to the 141-node power network under the same conditions detailed in Section \ref{sec:bin_2stage}. This creates a large-scale LP with 1,420,141 variables and 4,230,001 constraints. The indicated values of $y^k$ are again all verified to be integer feasible. These results are provided in Table \ref{tab:141_data} in the Appendix, and Figure \ref{fig:141node_pareto_cont} shows the Pareto pairs obtained from both formulations. In this case, the Pareto solutions are very similar but differences are more perceptible but rather small. Specifically, the identified values of $y^k$ differ by 0.49\% on average relative to the optimal values obtained via the mixed-integer formulation (i.e., 49 differences in the set of $y^k$ values out of the 10,000). The continuous Pareto pairs only required 13.41s on average to solve (in contrast to the mixed-integer results for which the majority were unable to converge within the 3,600s time limit).

Our results highlight that the continuous optimal design formulation provides high quality solutions with significant reductions in computational time. As mentioned in Section \ref{sec:cont_formulation}, the high quality of the solutions can likely be attributed to the degenerate nature of the SF index.

\begin{figure}[!htb]
		\centering
		\includegraphics[width=0.6\textwidth]{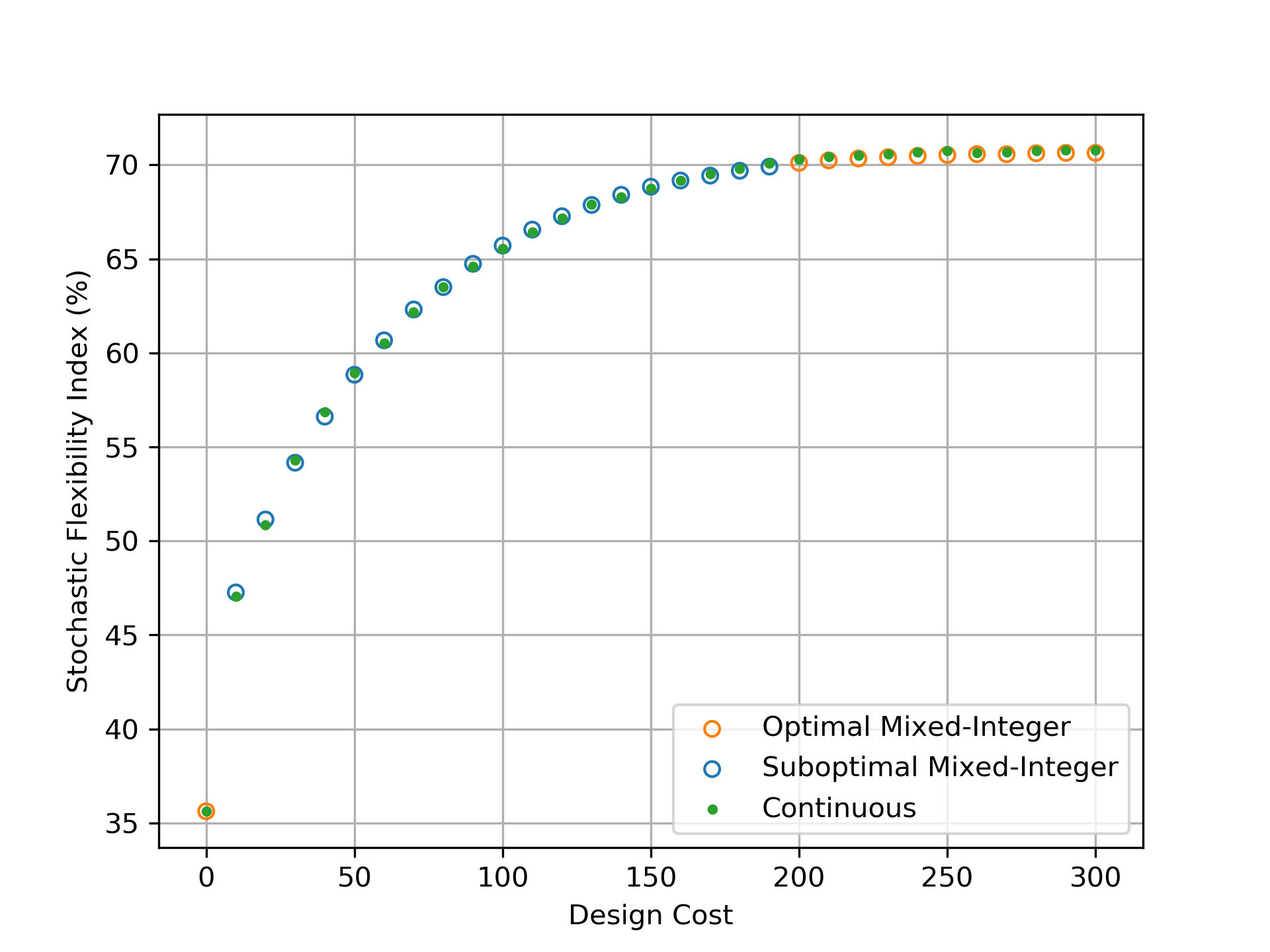}
	\caption{Pareto sets for 141-node power network  obtained with the mixed-integer and continuous formulations.}
	\label{fig:141node_pareto_cont}
\end{figure}

\FloatBarrier

\section{Conclusions and Future Work}
We have presented two stochastic optimal design formulations that minimize design cost and maximize or enforce design flexibility using the stochastic flexibility index. In particular, we are able to reproduce the results of joint chance constraint based formulations using a continuous conflict resolution program at a significantly reduced computational cost. In future work, we are interested in including recourse/operational costs and in gaining more theoretical insight as to why the continuous formulations provide high quality solutions and also to derive more effective rounding strategies.  Moreover, the continuous design formulation is structured an amenable to parallel decomposition. Consequently, we will explore the use of such techniques to enable the solution of large problems. 

\section*{Acknowledgments}
This work was supported by the U.S. Department of Energy under grant DE-SC0014114. We thank James Luedtke and Yue Shao for technical discussions.

\appendix 

\section{Supplemental Data}
The numerical results from Section \ref{sec:cases} are provided below in Tables \ref{tab:3d_data1}, \ref{tab:3d_data2}, \ref{tab:ieee14_data}, and \ref{tab:141_data}. Tables \ref{tab:3d_data1} and \ref{tab:3d_data2} feature the Pareto set data for the three-node network obtained with Problems \eqref{eq:2_stage_dist} and \eqref{eq:2_stage_mod_dist}. Table \ref{tab:ieee14_data} features the Pareto set data for the IEEE-14 node power network obtained with Problems \eqref{eq:2_stage_dist} and \eqref{eq:2_stage_mod_dist}. Table \ref{tab:141_data} features the Pareto set data for the 141-node power network with obtained with Problems \eqref{eq:2_stage_dist} and \eqref{eq:2_stage_mod_dist}.

\begin{table}[!htb]
	\caption{The first 22 results obtained for the three-node network using the mixed-integer and continuous formulations.}
	\begin{center}
		\begin{tabular}{|c|cccccc|}
			\hline	
			 \multirow{ 2}{*}{$\epsilon_c$} & Design & \multirow{ 2}{*}{$SF_K$ (\%)} & \multirow{ 2}{*}{$\bar{SF}_K$ (\%)} & Mixed-Integer & Continuous & Differences\\
			 & Cost &  &  & Time (s) & Time (s) & in $y^k$ (\%)\\ \hline \hline
			 0     & 0     & 96.7  & 96.7  & 0.0155 & 0.0156 & 0 \\
			 0.25  & 0.25  & 96.8  & 96.8  & 0.0170 & 0.0156 & 0 \\
			 0.5   & 0.5   & 97.0  & 97.0  & 0.0171 & 0.0000 & 0 \\
			 0.75  & 0.75  & 97.0  & 97.0  & 0.0179 & 0.0040 & 0 \\
			 1     & 1     & 97.2  & 97.2  & 0.0189 & 0.0157 & 0 \\
			 1.25  & 1.25  & 97.8  & 97.8  & 0.0190 & 0.0156 & 0 \\
			 1.5   & 1.5   & 97.9  & 97.9  & 0.0173 & 0.0000 & 0 \\
			 1.75  & 1.75  & 98.3  & 98.3  & 0.0195 & 0.0178 & 0 \\
			 2     & 2     & 98.3  & 98.3  & 0.0187 & 0.0040 & 0 \\
			 2.25  & 2.25  & 98.5  & 98.5  & 0.0168 & 0.0156 & 0 \\
			 2.5   & 2.5   & 98.6  & 98.6  & 0.0191 & 0.0000 & 0 \\
			 2.75  & 2.75  & 98.8  & 98.8  & 0.0175 & 0.0050 & 0 \\
			 3     & 3     & 98.9  & 98.9  & 0.0208 & 0.0000 & 0 \\
			 3.25  & 3.25  & 98.9  & 98.9  & 0.0203 & 0.0000 & 0 \\
			 3.5   & 3.5   & 99.0  & 99.0  & 0.0193 & 0.0000 & 0 \\
			 3.75  & 3.75  & 99.0  & 99.0  & 0.0180 & 0.0000 & 0 \\
			 4     & 4     & 99.1  & 99.1  & 0.0172 & 0.0000 & 0 \\
			 4.25  & 4.25  & 99.1  & 99.1  & 0.0200 & 0.0000 & 0 \\
			 4.5   & 4.5   & 99.3  & 99.3  & 0.0186 & 0.0156 & 0 \\
			 4.75  & 4.75  & 99.4  & 99.4  & 0.0192 & 0.0198 & 0 \\
			 5     & 5     & 99.5  & 99.5  & 0.0186 & 0.0050 & 0 \\
			 5.25  & 5.25  & 99.5  & 99.5  & 0.0182 & 0.0000 & 0 \\
			  \hline
		\end{tabular}
	\end{center}
	\label{tab:3d_data1}
\end{table}

\begin{table}[!htb]
	\caption{The last 21 results obtained for the three-node network using the mixed-integer and continuous formulations.}
	\begin{center}
		\begin{tabular}{|c|cccccc|}
			\hline	
			\multirow{ 2}{*}{$\epsilon_c$} & Design & \multirow{ 2}{*}{$SF_K$ (\%)} & \multirow{ 2}{*}{$\bar{SF}_K$ (\%)} & Mixed-Integer & Continuous & Differences\\
			& Cost &  &  & Time (s) & Time (s) & in $y^k$ (\%)\\ \hline \hline
			5.5   & 5.5   & 99.5  & 99.5  & 0.0183 & 0.0000 & 0 \\
			5.75  & 5.75  & 99.6  & 99.6  & 0.0193 & 0.0157 & 0 \\
			6     & 6     & 99.7  & 99.7  & 0.0177 & 0.0000 & 0 \\
			6.25  & 6.25  & 99.7  & 99.7  & 0.0187 & 0.0000 & 0 \\
			6.5   & 6.5   & 99.7  & 99.7  & 0.0202 & 0.0000 & 0 \\
			6.75  & 6.75  & 99.7  & 99.7  & 0.0179 & 0.0050 & 0 \\
			7     & 7     & 99.7  & 99.7  & 0.0185 & 0.0015 & 0 \\
			7.25  & 7.25  & 99.7  & 99.7  & 0.0207 & 0.0040 & 0 \\
			7.5   & 7.5   & 99.7  & 99.7  & 0.0196 & 0.0000 & 0 \\
			7.75  & 7.75  & 99.7  & 99.7  & 0.0182 & 0.0000 & 0 \\
			8     & 8     & 99.7  & 99.7  & 0.0183 & 0.0158 & 0 \\
			8.25  & 8.25  & 99.8  & 99.8  & 0.0195 & 0.0000 & 0 \\
			8.5   & 8.5   & 99.8  & 99.8  & 0.0192 & 0.0000 & 0 \\
			8.75  & 8.75  & 99.8  & 99.8  & 0.0192 & 0.0000 & 0 \\
			9     & 9     & 99.8  & 99.8  & 0.0207 & 0.0000 & 0 \\
			9.25  & 9.25  & 99.8  & 99.8  & 0.0192 & 0.0156 & 0 \\
			9.5   & 9.5   & 99.8  & 99.8  & 0.0182 & 0.0000 & 0 \\
			9.75  & 9.75  & 99.9  & 99.9  & 0.0180 & 0.0000 & 0 \\
			10    & 10    & 99.9  & 99.9  & 0.0202 & 0.0000 & 0 \\
			10.25 & 10.25 & 99.9  & 99.9  & 0.0201 & 0.0040 & 0 \\
			10.5  & 10.5  & 100   & 100   & 0.0195 & 0.0000 & 0 \\ \hline
		\end{tabular}
	\end{center}
	\label{tab:3d_data2}
\end{table}

\begin{table}[!htb]
	\caption{The results obtained for the IEEE-14 power network using the mixed-integer and continuous formulations.}
	\begin{center}
		\begin{tabular}{|c|cccccc|}
			\hline	
			\multirow{ 2}{*}{$\epsilon_c$} & Design & \multirow{ 2}{*}{$SF_K$ (\%)} & \multirow{ 2}{*}{$\bar{SF}_K$ (\%)} & Mixed-Integer & Continuous & Differences\\
			& Cost &  &  & Time (s) & Time (s) & in $y^k$ (\%)\\ \hline \hline
            0     & 0     & 88.10  & 88.10  & 16.20639 & 0.788219 & 0 \\
            2.5   & 2.5   & 88.30  & 88.30  & 47.31307 & 1.19029 & 0 \\
            5     & 5     & 88.95 & 88.95 & 36.28867 & 1.56477 & 0 \\
            7.5   & 7.5   & 89.20  & 89.20  & 51.17903 & 2.046263 & 0 \\
            10    & 10    & 89.50 & 89.50  & 54.72244 & 1.790163 & 0 \\
            12.5  & 12.5  & 89.80  & 89.80  & 84.38467 & 1.499423 & 0 \\
            15    & 15    & 90.15 & 90.15 & 94.45383 & 1.705751 & 0 \\
            17.5  & 17.5  & 90.30  & 90.30  & 71.09405 & 1.671522 & 0 \\
            20    & 20    & 90.55 & 90.55 & 49.37532 & 1.978597 & 0 \\
            22.5  & 22.5  & 90.70  & 90.70  & 72.26784 & 1.67794 & 0 \\
            25    & 25    & 90.80  & 90.80  & 36.79027 & 1.598156 & 0 \\
            27.5  & 27.5  & 90.95 & 90.95 & 34.87547 & 1.908839 & 0 \\
            30    & 30    & 91.05 & 91.05 & 45.3939 & 1.970099 & 0 \\
            32.5  & 32.5  & 91.20  & 91.20  & 57.60199 & 1.913504 & 0 \\
            35    & 35    & 91.20  & 91.20  & 50.41451 & 1.501976 & 0 \\
            37.5  & 37.5  & 91.20  & 91.20  & 25.38126 & 1.753275 & 0 \\
            40    & 40    & 91.25 & 91.25 & 20.42641 & 2.240765 & 0 \\
            42.5  & 42.5  & 91.30  & 91.30  & 19.73717 & 1.493643 & 0 \\
            45    & 45    & 91.30  & 91.30  & 22.286 & 1.579459 & 0 \\
            47.5  & 47.5  & 91.35 & 91.35 & 22.73862 & 1.853226 & 0 \\
            50    & 50    & 91.40  & 91.40  & 31.93482 & 2.758479 & 0 \\
            52.5  & 52.5  & 91.40  & 91.40  & 29.63871 & 1.616043 & 0 \\
            55    & 55    & 91.45 & 91.45 & 26.9819 & 1.736844 & 0 \\
            57.5  & 57.5  & 91.50  & 91.50  & 34.11424 & 1.60253 & 0 \\
            60    & 60    & 91.50  & 91.50  & 19.58831 & 2.56982 & 0 \\
            62.5  & 62.5  & 91.50  & 91.50  & 16.80612 & 1.569985 & 0 \\
            65    & 65    & 91.50  & 91.55 & 20.63473 & 1.492194 & 0.05 \\ \hline
		\end{tabular}
	\end{center}
	\label{tab:ieee14_data}
\end{table}

\begin{table}[!htb]
	\caption{The results obtained for the 141-node power network using the mixed-integer and continuous formulations.}
	\begin{center}
		\begin{tabular}{|c|ccccccc|}
			\hline	
			\multirow{ 2}{*}{$\epsilon_c$} & Design & \multirow{ 2}{*}{$SF_K$ (\%)} & \multirow{ 2}{*}{$\bar{SF}_K$ (\%)} & Mixed-Integer & Continuous & Differences & Optimal\\
			& Cost &  &  & Time (s) & Time (s) & in $y^k$ (\%) & Mixed-Integer\\ \hline \hline
            0     & 0     & 35.62 & 35.62 & 7.09  & 12.52 & 0     & Yes \\
            10    & 10    & 47.26 & 47.07 & 3605.73 & 13.66 & 0.95  & No \\
            20    & 20    & 51.14 & 50.84 & 3608.59 & 14.18 & 2.02  & No \\
            30    & 30    & 54.16 & 54.29 & 3614.13 & 13.77 & 1.95  & No \\
            40    & 40    & 56.61 & 56.86 & 3617.04 & 11.57 & 2.17  & No \\
            50    & 50    & 58.84 & 58.93 & 3622.68 & 13.46 & 1.79  & No \\
            60    & 60    & 60.67 & 60.53 & 3625.90 & 15.61 & 2.24  & No \\
            70    & 70    & 62.31 & 62.19 & 3625.88 & 15.21 & 2.12  & No \\
            80    & 80    & 63.50  & 63.50  & 3641.28 & 12.29 & 1.68  & No \\
            90    & 90    & 64.74 & 64.60  & 3618.84 & 15.34 & 1.42  & No \\
            100   & 100   & 65.71 & 65.54 & 3641.04 & 15.49 & 1.97  & No \\
            110   & 110   & 66.56 & 66.42 & 3624.48 & 12.30 & 2.1   & No \\
            120   & 120   & 67.27 & 67.18 & 3620.28 & 13.01 & 1.79  & No \\
            130   & 130   & 67.87 & 67.89 & 3628.07 & 15.14 & 1.3   & No \\
            140   & 140   & 68.41 & 68.29 & 3619.19 & 12.63 & 1.24  & No \\
            150   & 150   & 68.84 & 68.76 & 3619.67 & 15.29 & 1.2   & No \\
            160   & 160   & 69.17 & 69.18 & 3635.72 & 12.64 & 1.13  & No \\
            170   & 170   & 69.43 & 69.52 & 3624.00 & 13.38 & 1.01  & No \\
            180   & 180   & 69.69 & 69.81 & 3623.32 & 14.87 & 0.76  & No \\
            190   & 190   & 69.91 & 70.10  & 3623.43 & 11.80 & 0.59  & No \\
            200   & 200   & 70.11 & 70.28 & 804.34 & 11.51 & 0.49  & Yes  \\
            210   & 210   & 70.24 & 70.43 & 193.87 & 13.10 & 0.47  & Yes  \\
            220   & 220   & 70.34 & 70.50  & 160.40 & 12.64 & 0.5   & Yes  \\
            230   & 230   & 70.41 & 70.59 & 165.22 & 12.67 & 0.5   & Yes  \\
            240   & 240   & 70.48 & 70.67 & 124.44 & 12.42 & 0.49  & Yes  \\
            250   & 250   & 70.53 & 70.76 & 111.72 & 12.45 & 0.47  & Yes  \\
            260   & 260   & 70.57 & 70.64 & 108.32 & 14.97 & 0.63  & Yes  \\
            270   & 270   & 70.57 & 70.69 & 103.03 & 13.21 & 0.64  & Yes  \\
            280   & 280   & 70.62 & 70.76 & 101.47 & 13.25 & 0.56  & Yes  \\
            290   & 290   & 70.64 & 70.78 & 101.27 & 13.55 & 0.54  & Yes  \\
            300   & 300   & 70.64 & 70.78 & 99.66 & 11.79 & 0.54  & Yes  \\ \hline
		\end{tabular}
	\end{center}
	\label{tab:141_data}
\end{table}

\bibliography{references}

\end{document}